\documentclass{article}
\usepackage{amsmath}
\usepackage[psamsfonts]{amssymb}
\usepackage{cmmib57}
\newcommand{\bPf}{\par\vspace*{-4pt}\indent{\sc Proof.}\enskip}
\newcommand{\ePf}{\medskip}
\def\QED{\hskip0.1em\hfill\null\ \null\nobreak\hfill\kern3pt\vbox{\hrule\hbox
   {\vrule\kern1pt\vbox{\kern1.7pt\hbox{$\scriptscriptstyle{QED}$}
    \kern0.2pt}\kern1pt\vrule}\hrule}}

\def\END{\hskip0.1em\hfill\null\ \null\nobreak\hfill\kern3pt\vbox{\hrule\hbox
   {\vrule\kern1pt\vbox{\kern1.7pt\hbox{$\,\,\,\vspace{5pt}$}
    \kern0.2pt}\kern1pt\vrule}\hrule}}
\newtheorem{theorem}{Theorem}
\newtheorem{lemma}{Lemma}
\newtheorem{corollary}{Corollary}
\newtheorem{proposition}{Proposition}
\newtheorem{remark}{Remark}
\newtheorem{definition}{Definition}
\newtheorem{example}{Example}
\newcommand{\bCd}{\bEq\begin{CD}}
\newcommand{\eCd}{\end{CD}\eEq}
\newcommand{\bcd}{\beq\begin{CD}}
\newcommand{\ecd}{\end{CD}\eeq}
\newcommand{\ben}{\begin{enumerate}}
\newcommand{\een}{\end{enumerate}}
\newcommand{\bEq}{\begin{eqnarray}}
\newcommand{\eEq}{\end{eqnarray}}
\newcommand{\beq}{\begin{eqnarray*}}
\newcommand{\eeq}{\end{eqnarray*}}
\newcommand{\bDf}{\begin{definition}\em}
\newcommand{\eDf}{\end{definition}}
\newcommand{\bLm}{\begin{lemma}}
\newcommand{\eLm}{\end{lemma}}
\newcommand{\bPr}{\begin{proposition}}
\newcommand{\ePr}{\end{proposition}}
\newcommand{\bTh}{\begin{theorem}}
\newcommand{\eTh}{\end{theorem}}
\newcommand{\bCr}{\begin{corollary}}
\newcommand{\eCr}{\end{corollary}}
\newcommand{\bRm}{\begin{remark}\em}
\newcommand{\eRm}{\end{remark}}
\newcommand{\bEx}{\begin{example}\em}
\newcommand{\eEx}{\end{example}}
\newcommand{\C}{\mathbb{C}}

\newcommand{\Z}{\mathbb{Z}}
\newcommand{\ie}{{\em i.e.} }
\newcommand{\eg}{{\em e.g.} }

\newcommand{\mto}{\mapsto}

\newcommand{\btht}{\boldsymbol{\tht}}

\newcommand{\cD}{\mathcal{D}}

\newcommand{\cI}{\mathcal{I}}

\newcommand{\bP}{\boldsymbol{P}}
\newcommand{\bQ}{\boldsymbol{Q}}

\newcommand{\bS}{\boldsymbol{S}}

\newcommand{\bX}{\boldsymbol{X}}

\newcommand{\sub}{\subset}

\newcommand{\wed}{\wedge}

\newcommand{\alp}{\alpha}
\newcommand{\bet}{\beta}
\newcommand{\gam}{\gamma}
\newcommand{\del}{\delta}
\newcommand{\eps}{\epsilon}
\newcommand{\zet}{\zeta}
\newcommand{\tht}{\theta}

\newcommand{\kap}{\kappa}

\newcommand{\sig}{\sigma}

\newcommand{\ome}{\omega}
\newcommand{\Gam}{\Gamma}

\newcommand{\Sig}{\Sigma}
\newcommand{\Ome}{\Omega}

\title{\large{{\bf  B\"acklund loop algebras for compact
 and non-compact nonlinear spin models
in $(2+1)$ dimensions\thanks{This note is part of a joint research work with E. Winterroth.}}}}
\author{{\normalsize M.
Palese
}\thanks{Partially
supported by
GNFM of INdAM, MURST (PRIN 2003) and University of Torino.}
\\{\footnotesize Department of Mathematics,
University of Torino}
\\{\footnotesize Via C. Alberto 10, 10123 Torino, Italy}\\ 
{\footnotesize e--mail: 
{\sc marcella.palese@unito.it
}}\\
}
\date{}
\begin{document}

\maketitle

\begin{abstract}
The B\"acklund problem is solved for both
the compact and noncompact versions of the Ishimori (2+1)-dimensional nonlinear spin model. 
In particular, a realization of the arising B\"acklund algebra
in the form of an infinite-dimensional loop Lie algebra of the Ka\v{c}--Moody type is provided.

\medskip

{\footnotesize \noindent \small{{\bf Key words}: 
 integrable systems, nonlinear spin models, prolongation algebras, 
 B\"acklund transformations, B\"acklund-Cartan connections.

}}

\end{abstract}

\section{Introduction}

In the study of nonlinear field equations, the 
B\"acklund structures method, also known as prolongation structures method (see \eg
\cite{Wa73,Wa75,Mo76,PRS79,DoFo83,To85,Fo90,Pa93,ALLPS94,ALLPS95,PLS00,PaWi02,PaWi03} for a review of the
procedure in some relevant applications concerning both the (1+1) and (2+1) dimensional cases)
constitutes
a systematic analytical procedure which enables one, in principle, to associate a linear 
problem with the equation under consideration. Within such a method, 
nonlinear (prolongation) B\"acklund algebras are related to integrable nonlinear 
field equations which can be expressed by means of closed differential ideals.
Such algebras arise via the introduction of an arbitrary number of 
B\"acklund forms containing new dependent variables (called 
pseudopotentials), and by requiring the algebraic equivalence between the 
generators of the prolonged ideal and its exterior differential, \ie by requiring an
integrability condition for the prolonged differential ideal.
It was pointed out that integrability properties of nonlinear field equations can be related with the
existence of B\"acklund symmetries and admissible B\"acklund transformations (see {\em e.g.}
\cite{AI79,AbSe81}). In fact, B\"acklund algebras appear in form of incomplete Lie algebra
structures (in  the sense that not all of the commutators are known) 
as necessary conditions for integrability of
a connection induced by a B\"acklund map \cite{PRS79}, and a realization of such algebras (particularly
in loop--algebras form) 
provides the solution of the so--called B\"acklund problem. 

In more than two independent variables, the extension of
the prolongation procedure is generally nontrivial and some aspects remain to be explored (see, for
example, \cite{Mo76,To85,PLS00,PaWi02}). Recently, relying on a suitable  geometric description
of B\"acklund transformations \cite{AI79}, some developments have been achieved on the theoretical side, which
provide some hints for a better  understanding of the whole matter \cite{PaWi03}. The study of higher
dimensional systems is a central theme in the theory of integrable systems; in \cite{PaWi02} the general
approach of generating new (2+1)--dimensional integrable systems from a given abstract algebraic structure
via the extension of a B\"acklund map was tackled by showing
how a connection can be induced by a B\"acklund map in the jet bundles framework and a
characterization of completely integrable systems in terms of B\"acklund 
structures has been provided. In \cite{PaWi03} we also proved that such a connection is 
an admissible Cartan connection for a suitable $\bar{K}$-structure, \ie a subbundle of a (higher order) frame
bundle.  

In this note, making use of an {\em ansatz} \cite{Mo76,PRS79,To85,PaWi02} 
which is in fact a suitable version 
of the structure equations for an admissible B\"acklund-Cartan connection \cite{PaWi03},
we associate an infinite-dimensional B\"acklund algebra structure
with the Ishimori $(2+1)$-dimensional 
spin model in both the compact \cite{Ish84} and the noncompact versions \cite{LMS92}. Homomorphisms with
quotient finite and infinite dimensional Lie algebras are worked out. In particular, in view of possible
applications - \eg to derive a whole family of (2+1)-dimensional nonlinear field equations
containing the original Ishimori spin models -, a realization of the arising B\"acklund algebra
in the form of an infinite-dimensional loop Lie algebra of the Ka\v{c}--Moody type is provided.

\section{Admissible B\"acklund transformations}
In the following we shall shortly recall few basics concepts and set the notation.
We shall assume the reader is familiar with the basic notions from 
the theory of bundles, jet prolongations, principal bundles 
and connections (for references and details see \eg \cite{Ko72,Sau89}).

Let $\pi:  U\to  X$, $\tau: Z\to  X$,
be two (vector) bundles 
with local fibered coordinates $(x^{\alp},u^{A})$ and $(x^{\alp},z^{i})$,
respectively, where $\alp=1,\ldots,m=\textstyle{dim} X$, $A=1,\ldots,n=\textstyle{dim} U -
\textstyle{dim} X$, $i=1,\ldots,N=\textstyle{dim}Z-\textstyle{dim} X$.  
A system of nonlinear field equations of order $k$ on $ U$ is geometrically described as an
exterior differential system $\nu$ on $J^{k} U$. The solutions of the field
equations are (local)  sections $\sig$
of $U\to X$ such that $(j^{k}\sig)^{*}\nu=0$. We shall also denote by
$J^{\infty}\nu$ ({\em resp.} $j^{\infty}\sig$) 
the infinite order jet prolongation of $\nu$ ({\em resp.} $\sig$). 

Let then $B$ be the infinite--order contact transformations group on 
$J^{\infty}U$.

\bDf\label{admissible}
The group of (infinitesimal) {\em B\"acklund transformations for the
system $\nu$} is the closed subgroup
${ \tilde K}$ of $ B$ which leaves invariant solution submanifolds 
of $J^{\infty}\nu$.
The group of (infinitesimal) {\em generalized B\"acklund transformations for the system
$\nu$} is the closed subgroup $K$ of
$B$ which leaves invariant $J^{\infty}\nu$ \cite{AI79}.\END\eDf

Let $\pi: U\to X$, $\tau:Z\to X$, be vector bundles as the above and 
$\pi^{1}: {J^{1} U} \to X$, $\tau^{1}: {J^{1}Z} \to X$, the first order jet
prolongations bundles, with local fibered coordinates
$(x^{\alp},u^A,u^{A}_{\alp})$,
$(x^{\alp},z^i,z^{i}_{\alp})$, respectively. Furthermore, let
$(dx^{\beta}$, 
$du^{A}$, $du^{A}_{\beta})$ and $(dx^{\beta},dz^{i}$, $dz^{i}_{\beta})$ be local
bases of $1$--forms on $J^{1} U$ and $J^{1}Z$,
respectively.

\bDf
We define a first order B\"acklund map to be 
the fibered morphism over $Z$:
$\phi: J^{1} U\times_{X} Z\to J^{1}Z: (x^{\alp},u^A,u^{A}_{\alp};z^{i})\mto
(x^{\alp},z^{i},z^{i}_{\alp})$,
with $z^{i}_{\alp}=\phi^{i}_{\alp}(x^{\bet},u^A,u^{A}_{\bet};z^{j})$.

The fibered morphism $\phi$ is said to be an {\em admissible} B\"acklund
transformation for the differential system $\nu$ if
$\phi^{i}_{\alpha}=\cD_{\alpha}\phi^{i}$ - where $\cD_{\alpha}$ 
is the formal derivative on $J^{1} U\times_{X}Z$
- and the integrability conditions coincide with the
exterior differential system $\nu$.
\END\eDf

\bRm 
By pull--back of the contact structure on $J^1 Z$, the B\"acklund morphism 
induces an horizontal distribution, the
{\em induced B\"acklund connection},
on the bundle $(J^{1}{ U}\times_{X}{Z},$ $J^{1}{ U},$
$\pi^{1*}_{0}(\eta))$ \cite{PRS79}.\END\eRm

\bDf
Let $\bar{K}$ be a normal subgroup
$\bar{K}\sub ({\tilde K} \cap K) \sub  B$ leaving invariant (the infinite order
prolongation of) $\nu$ and its solutions. We call $\bar{K}$ the group of (infinitesimal) {\em generalized
admissible B\"acklund transformations} for the system $\nu$.
\END\eDf

\bTh\label{Th1}
The following statements are equivalent \cite{PaWi02}.
\begin{enumerate}
\item $\phi$ is an admissible B\"acklund transformation for the
differential system $\nu$. 
\item The induced B\"acklund connection is
$\bar{K}$--invariant.
\end{enumerate}
\eTh

\section{B\"acklund algebras for Ishimori models}

Let us now consider the nonlinear spin model originally introduced by 
Ishimori \cite{Ish84} to extend in (2+1) dimensions the continuous isotropic Heisemberg spin model in 
(1+1) dimensions (for the study of the integrability properties of 
the latter see \eg \cite{Lak77,Tak77,ALLPS94}). 
 A Lax pair was provided by Ishimori and multivortex solutions were derived by means of
the Hirota method thus showing that the vortex dynamics is integrable.
In \cite{LMS92} was shown the gauge equivalence between a noncompact version of the Ishimori spin model and
the Davey--Stewartson equation. 

The classical continuous isotropic Heisenberg spin chain can be generalized to the following set of
coupled nonlinear
$(2+1)$-dimensional partial differential equations:
\bEq
\Sig {\bS}  _{t}={\bS}   \times 
({\bS}  _{xx}+\eps^{2}{\bS}  _{yy})+\phi_{y}{\bS}  _{x}+\phi_{x}{\bS}  _{y}\,,\label{eq1} \\
\phi_{xx}-\eps^{2}\phi_{yy}=-2\eps^{2}\Sig {\bS}\,\cdot 
\,({\bS}  _{x}\times {\bS}  _{y})\,, \label{eq2} 
\eEq
where ${\bS}  ={\bS}  (x,y,t)$ is a classical spin field vector, 
$\Sig$ is a diagonal matrix $diag$ $(1,1,\kap^2)$, $\kap^2=\pm1$, 
subscripts denote partial derivatives, and the symbol $\times$ stands for 
the usual vector product. The spin field components $S_j$, $j=1,2,3$ are assumed to 
satisfy the constraint
\bEq\label{constraint}
\Sig{\bS}\, \cdot\,{\bS}  =\kap^2\,
\eEq 
where $\kap^2=\pm1$ refers to the compact \cite{Ish84} and the noncompact \cite{LMS92} case, 
respectively. In other words, for $\kap^2=1$ the quantities $S_j$ belong 
to the unitary sphere $SU(2)/U(1)$, while for $\kap^2=-1$, the $S_j$'s
range over a shield of the two-fold hyperboloid $SU(1,1)/U(1)$.

Let us introduce on a (vector) bundle $J^{1}U$ with local coordinates
$(x,y,t;$ ${\bS}  ,{\bP} ,{\bQ};$ $\phi,\alp, \mu)$ \footnote{For convenience
we made the following change of notation: $x^{1}=x$, $x^{2}=y$, $x^{3}=t$, $u^{1}={\bS}  $,
$u^{1}_{1}={\bP} $, $u^{1}_{2}={\bQ}$, $u^{2}=\phi$, $u^{2}_{1}=\alp$, $u^{2}_{2}=\mu$; furthermore, in the
sequel, $z=\xi$.} the {\em closed} differential ideal defined by the set  of (vector-valued)
$3$-forms:
\bEq
& &{\btht}_{1}=(d {\bS}-{\bP}  dx)\wed dy\wed dt\,,  \label{ideal1}
\\ 
& &{\btht}_{2}=(d {\bS}-{\bQ} dy)\wed dx\wed dt\,, \label{ideal2}
\\ 
& & {\btht}_{3}=d\Sig {\bS} \wed dy\wed dx
+{\bS}  \times (d {\bP} \wed dy\wed 
dt-\eps^{2}d{\bQ} \wed dx\wed 
dt)+ \nonumber \\ 
& &+(\mu{\bP} +\alp{\bQ} )dx\wed dy\wed dt \,,
\label{ideal3}
\eEq
by the scalar $3$-forms:
\bEq
& &\gam_{1}=(d \phi-\alp dx)\wed dy\wed dt\,,  \label{ideal4}
\\
& &\gam_{2}=(d\phi-\mu dy)\wed dx\wed dt\,, \label{ideal5}
\\ 
& & \gam_{3}=d\alp \wed dy\wed dt +\eps^{2}d\mu \wed dx\wed 
dt+  \nonumber \\
 & &+2\eps^{2}\Sig {\bS}\, \cdot\, ({\bP} \times {\bQ})dx\wed dy
\wed dt \,, \label{ideal6}
\eEq
and by
\bEq
& &\bet_{1}=(d\Sig{\bS}\,\cdot\, {\bP} + \Sig{\bS}\,\cdot\, d{\bP}
)\wed dy\wed  dt\,, \label{ideal7} \\ 
& &\bet_{2}=(d\Sig{\bS}\,\cdot\, {\bQ} + \Sig{\bS}\,\cdot\, d{\bQ} )
\wed dx\wed dt\,, 
\label{ideal8}
\eEq
where $\wed$ stands for the exterior product of forms.

It is easy to verify the following.
\bPr
On every integral submanifold defined by (local) sections ${\bS}   = {\bS}  (x,y,t)$, 
${\bP} ={\bS}  _{x}$, ${\bQ}={\bS} _{y}$,
$\phi=\phi(x,y,t)$, 
$\alp=\phi_{x}$, $\mu=\phi_{y}$
- for which $dx\wedge dy\wedge dt\neq 0$, because of transversality of fibers -, the ideal
(\ref{ideal1})--(\ref{ideal8}) is equivalent to Eqs. (\ref{eq1})--(\ref{eq2}) with 
the constraint (\ref{constraint}). 
\ePr

Now let us consider the (vector) bundle $Z$ (fibered over the same basis as the above) with local coordinates
$(x,y,t; \xi^{m})$ and, defined on the fibered product $J^{1} U\times_{X}Z$, the $2$--forms:
\bEq\label{prolform}
\Ome^{k}=H^{k}dx\wed dy+F^{k}dy\wed dt+G^{k}dx\wed 
dt+ \nonumber \\
+(A^{k}_{m}dx+B^{k}_{m}dy+\del^{k}_{m}dt)\wed d\xi^{m}\,,
\eEq
where $\xi=\{\xi^{m}\}$, $k,\,m\,=\,1,2,\ldots,
{\rm N}$ (N at this stage is arbitrary), and $H^{k}$, $F^{k}$ and $G^{k}$ are, respectively, the
pseudopotential and functions on $J^{1} U\times_{X}Z$ to be determined. Furthermore, 
the quantities $A^{k}_{m}$
and $B^{k}_{m}$ denote the elements of two $N \times N$ 
constant regular matrices, $\del^{k}_{m}$ is the Kronecker symbol and
the summation convention over repeated indices is understood.

\bDf
The forms $\Omega^k$ are called the {\em B\"acklund forms} associated 
with Eqs. (\ref{eq1})--(\ref{constraint}).
Let $\cI$ be the {\em prolonged} ideal generated by the forms 
${\btht}_{i},\gam_{i},\bet_{l}$, $i=1,2,3$, $l=1,2$ and
$\Omega^{k}$, $k=1,\ldots,N$.  We say that $\cI$ is closed if 
$d\Omega^{k}\in \cI$
$({\btht}_{i},\gam_{i},\bet_{l},\Omega^{k})$.
\END\eDf 

\bRm
The reason for we call such prolongation forms `B\"acklund' is that Eq. (\ref{prolform})
is in fact a specific form of the generalized {\em ansatz} introduced in \cite{PaWi02}, which is assumed to be
induced by an admissible B\"acklund transformation. The induced B\"acklund-Cartan connection on
$(J^{1}{ U}\times_{X}{Z},$
$J^{1}{ U},$
$\pi^{1*}_{0}(\eta))$ is here given by
\bEq\label{Baec}
\ome^{m}=\Gam^{m}_{1}dx+\Gam^{m}_{2}dy+\Gam^{m}_{3}dt+d\xi^{m} \,,
\eEq
and the closed (matrix-valued) form $\theta$ of \cite{PaWi02} is here given by
$\theta^{k}_m =A^{k}_{m}dx+B^{k}_{m}dy+\del^{k}_{m}dt$. With this notation in mind then
$H^{k}= \Gam^{m}_{1}B^{k}_{m}-\Gam^{m}_{2}A^{k}_{m}$, $F^{k}= \Gam^{m}_{2}\del^{k}_{m}-\Gam^{m}_{3}B^{k}_{m}$,
$G^{k}=\Gam^{m}_{1}\del^{k}_{m}-\Gam^{m}_{3}A^{k}_{m}$.
The closure condition for the ideal
$\cI$ is an equivalent condition to Theorem
\ref{Th1}, where $\bar{K}$ is here to be determined by solving the so-called B\"acklund problem. The B\"acklund
map, and thus the induced Cartan connection, are completely determined by the functions
$H^{k}$,
$F^{k}$ and
$G^{k}$ and the matrices $\{A^{k}_{m}\}$
and $\{B^{k}_{m}\}$.\END
\eRm

The dependence of the functions $H^{k}, F^{k}, G^{k}$ 
on $({\bS}  ,{\bS}  _{x},{\bS}  _{y};\phi,\phi_{x}, 
\phi_{y};\xi^{m})$ shall be determined by requiring the B\"acklund transformation defined by Eq. (\ref{Baec})
to be admissible for the Ishimori systems, \ie the closure of the prolonged ideal to be satisfied.

In the following $[G, H]^{k} = G^{j}H^{k}_{\xi^{j}} - H^{j}G^{k}_{\xi^{j}}$ (Lie bracket),
$H^{k}_{\xi^{j}} = \frac{\partial H^{k}}{\partial \xi^{j}}$, 
$H^{k}_{u} = \frac{\partial H^{k}}{\partial u}$, and so on.
Furthermore, we shall omit the indices $k$, $m$ for simplicity. The 
Lie brackets must be understood as commutators then.

\bLm\label{integrability}
The closure condition for $\cI$ yields:
\bEq
& &H({\bS}  ;\xi)={\bX}(\xi)\,\cdot\,{\bS}  +Y(\xi)\,, \\
& &F=-({\bX}\times{\bS})\,\cdot\,{\bS}  _{x}+\hat{F}({\bS}  ;;\xi)\,,\\
& &G=\eps^{2}({\bX}\times{\bS})\,\cdot\,{\bS}  _{y}+\hat{G}({\bS}  ;\xi)\,\\
& &\hat{F}_{{\bS}  }\,\cdot\,{\bS}  _{x}-\hat{G}_{{\bS}  }\,
\cdot\,{\bS}  _{y}+[{\bX}(\xi)\,\cdot\,{\bS}   +Y(\xi), 
\eps^{2}({\bX}\times{\bS})\,\cdot\,{\bS}  _{y}]+ \nonumber\\
& &+[{\bX}(\xi)\,\cdot\,{\bS}   +Y(\xi),\bar{G}({\bS}  ;\phi,\phi_{y};\xi)]=0\,.
\eEq
\eLm

\bPf
The closure condition is equivalent to the following constraints:
\bEq\label{intermezzo}
& & H_{\phi}=H_{\phi_{x}}=H_{\phi_{y}}=G_{\phi_{x}}=F_{\phi_{y}}=0\,,\\ 
& & H_{{\bS}  _{x}}=G_{{\bS}  _{x}}=0\,, \quad 
F_{{\bS}  _{x}}=-(\Sig H_{{\bS}  })\times\bS\,,\\
& & H_{{\bS}  _{y}}=F_{{\bS}  _{y}}=0\,, \quad 
G_{{\bS}  _{y}}=\eps^{2}(\Sig H_{{\bS}  })\times{\bS}\,,\\  
& &G_{\phi_{y}}-\eps^{2}F_{\phi_{x}}=0\,,
\eEq
\bEq\label{relations}
& &[H,G]+H_{{\bS}  }\,\cdot\,(\phi_{y}{\bS}  _{x}+\phi_{x}{\bS}  _{y})+F_{{\bS}  }\,\cdot\,{\bS}  _{x}-
G_{{\bS}  }\,\cdot\,{\bS}  _{y}
+ \nonumber \\
& &+ F_{\phi}\,\cdot\,\phi_{x}-G_{\phi}\,\cdot\,\phi_{y}-2
\eps^{2}F_{\phi_{x}}\Sig{\bS}\,\cdot\,({\bS}  _{x}\times{\bS}_{y}) =0\,.
\eEq
\bEq
&\hphantom{=}& F_{\xi} - G_{\xi}A
- H_{\xi}B= 0 \,,\label{spectral1}\\
&\hphantom{=}& F_{\xi}(B^{-1})G 
- G_{\xi}(B^{-1})F = 0\,,\label{spectral2}
\eEq 
\bEq 
\left[A,B\right] = 0\,.\label{spectral3}
\eEq

From equations (\ref{intermezzo}) we get
\bEq
& & H=H({\bS}  ;\xi)\,, \quad F=F({\bS}  ,{\bS}  _{x};\phi,\phi_{x};\xi)\,, \quad 
G=G({\bS}  ,{\bS}  _{y};\phi,\phi_{y};\xi)\,\\
& & F_{{\bS}  _{x}}=-\Sig H_{{\bS}  }\times{\bS}\,,\quad 
G_{{\bS}  _{y}}=\eps^{2}\Sig H_{{\bS}  }\times{\bS}\,,\quad 
G_{\phi_{y}}=\eps^{2}F_{\phi_{x}}\,.
\eEq
From the last three equations above we further get
\bEq\label{xx}
& & F=-(\Sig H_{{\bS}  }\times{\bS})\,\cdot\,{\bS}  _{x}+\bar{F}({\bS}  ;\phi,\phi_{x};\xi)\,,\\
& & 
G=\eps^{2}(\Sig H_{{\bS}  }\times{\bS})\,\cdot\,{\bS}  _{y}+\bar{G}({\bS}  ;\phi,\phi_{y};\xi)\,;\\
& & \bar{G}_{\phi_{y}}=\eps^{2}\bar{F}_{\phi_{x}}\,;
\eEq
The first two imply 
\bEq
H_{{S_{i}}_{{\bS}  }}=0\,,
\eEq
thus
\bEq\label{fund}
H({\bS}  ;\xi)={\bX}(\xi)\,\cdot\,{\bS}   +Y(\xi)\,,
\eEq
where ${\bX}=(X_{1},X_{2},X_{3})$.

Furthermore, we have $\bar{F}_{\phi_{x}}=0$, so that 
$\bar{F}=\bar{F}({\bS}  ;\phi;\xi)$.
Thus Eq. (\ref{relations}) become
\bEq
& & {\bX}\,\cdot\,{\bS}  _{x}\phi_{y}+{\bX}\,\cdot\,{\bS}  _{y}\phi_{x}+
\bar{F}_{{\bS}  }\,\cdot\,{\bS}  _{x}-\bar{G}_{{\bS}  }\,
\cdot\,{\bS}  _{y}+\bar{F}_{\phi}\phi_{x}-\bar{G}_{\phi}\phi_{y}+ \nonumber \\
& & 
+[{\bX}\,\cdot\,{\bS}   +Y, 
\eps^{2}({\bX}\times {\bS})\,\cdot\,{\bS}  _{y}+\bar{G}({\bS}  ;\phi,\phi_{y};\xi)]=0\,.
\eEq
From the equation above we further infer $\bar{F}_{\phi}=0$ and from 
Eqs. (\ref{xx}), $\bar{G}_{\phi_{y}}=0$, which in turns 
implies that $\bar{G}_{\phi}=0$.

So we finally get the assertion.
\ePf
\QED

In the following we provide the B\"acklund algebra $\bar{\mathfrak{k}}$ associated 
with the compact and non-compact Ishimori spin models.

\bPr
The following 
incomplete Lie algebra structure $\bar{\mathfrak{k}}$ is associated with Eqs.
(\ref{eq1})--(\ref{eq2}) and constraint
(\ref{constraint}):
\bEq
& &[Y,X_{1}]=[Y,X_{2}]=[Y,X_{3}]=0 \,, \\
& &[X_{1},[X_{1},X_{3}]]-[X_{2},[X_{2},X_{3}]]=0 \,, \\
& &k^{2}[X_{1},[X_{1},X_{2}]]+[X_{3},[X_{2},X_{3}]]=0 \,, \\
& &k^{2}[X_{2},[X_{1},X_{2}]]-[X_{3},[X_{1},X_{3}]]=0 \,, \\
& &[X_{1},[X_{2},X_{3}]]+[X_{2},[X_{1},X_{3}]]=0 \,, \\
& &[X_{2},[X_{1},X_{3}]]+[X_{3},[X_{1},X_{2}]]=0 \,, \\
& &k^{2}[X_{2},[X_{1},X_{3}]]+[Z,Y]=0 \,.
\eEq
\ePr

\bPf
From the above Lemma we get
\bEq\label{prealgebra}
& &[H({\bS}  ;\xi),\hat{G}({\bS};\xi)]=0\,,\\
& &\hat{F}_{{\bS}  }=0\,,
\eEq
which implies $\hat{F}=
K(\xi)$.
Thus we finally have to work out the equation:
\bEq
\hat{G}_{\bS} \,\cdot\,{\bS} _{y}=\eps^{2}[{\bS}\,\cdot\,{\bS}   
+Y,({\bX}\times{\bS})\,\cdot\,{\bS} _{y}]\,.
\eEq

By using the constraint (\ref{constraint}) we easily get 
\bEq
& &[Y,X_{1}]=[Y,X_{2}]=[Y,X_{3}]=0\,, \\
& &\hat{G}=-S_{1}[X_{2},X_{3}]+S_{2}[X_{1},X_{3}]-k^{2}S_{3}[X_{1},X_{2}] + 
Z(\xi)\,.
\eEq

Finally, the B\"acklund algebra $\bar{\mathfrak{k}}$ is obtained inserting the last result in Eq.
(\ref{prealgebra}) and using Eqs.
(\ref{xx}) and (\ref{fund}).
\QED\ePf

\bRm
The compatibility relations (\ref{spectral1}), (\ref{spectral2}) and (\ref{spectral3})
have to be satisfied by $K(\xi)$. It is easy to realize then that commutation relations for $K(\xi)$ also
hold true.\END
\eRm

It is possible to provide suitable realizations (\ie homomorphisms) of the algebra $\bar{\mathfrak{k}}$ with
finite or infinite dimensional quotient Lie algebras.

\bPr
The B\"acklund algebra $\bar{\mathfrak{k}}$
is homomorphic to the $\mathfrak{sl}(2,\C)$ algebra

\bEq
[X_1, X_2]=2i\zet \kappa^2 X_3, [X_1, X_3]=-2i\zet X_2, [X_2, X_3] = 
2i\zet X_1 \,
\eEq
where $\zet$ is a free parameter. 
\ePr

\bPf 
Let $[X_1, X_2]=X_4$, $[X_1, X_3]=X_5$, $[X_2, X_3]=X_6$ and assume that 
$X_1$, $X_2$, and $X_3$ are independent and $X_4$, $X_5$, 
$X_6$, $Y$, $Z$ and $K$ are a linear combination of the former operators. In such 
a way one finds $X_4 = 2i\zet \kappa^2 X_3$,\quad $X_5 = -2i\zet X_2$,
$X_6 = 2i\zet X_1$, and $Y=Z=K=0$. \QED\ePf

\bPr
Assume $Y=Z=K=0$. Then the B\"acklund algebra $\bar{\mathfrak{k}}$
is homomorphic to the infinite 
dimensional Lie algebra of the Ka\v{c}-Moody type:
\bEq\label{KM}
[T_i^{(m)},T_j^{(n)}] = i\epsilon^{k}_{ij} T_k^{(m + n)}\,,
\eEq
$\epsilon^{k}_{ij}$ being the Ricci tensor.
\ePr

\bPf
In fact, let us suppose that
the commutators $[X_2, X_6]$, $[X_3, X_5]$, $[X_3, X_6]$, $[X_4, X_5]$, 
$[X_4, X_6]$ and $[X_5, X_6]$ are unknown. Hence, it is easy to see that a realization of the 
incomplete Lie algebra $\bar{\mathfrak{k}}$ is provided by
\bEq
& & X_1 =\kappa T_1^{(1)}, X_2 = \kappa T_2^{(1)}, X_3 = T_3^{(1)} \,, \\
& & X_4 = i\kappa^2 T_3^{(2)}, X_5 = -i\kappa T_2^{(2)}, X_6 = i\kappa T_1^{(2)}\,, 
\eEq
where the vector fields $T_i^{(m)}$ $(i=1,2,3; m \in \Z)$ 
satisfy the commutation relations (\ref{KM}).
\QED\ePf

\bRm
It is easy to see that, in accordance with Theorem \ref{Th1}, the connection forms (\ref{Baec}) are
$\bar{\mathfrak{k}}$-invariant; in fact we have
$\Gam^{m}_i=\ome^{l}_i X^{m}_{l}$, with $i=1,2,3$ and $l=1,\ldots, M$; $M$ equals $\textstyle{dim}
\mathfrak{k}$ or the dimension of the quotient Lie algebra taken under consideration and $\ome{^l}_i$ are
admissible Cartan connection forms
\cite{PaWi03}.
\END\eRm

\subsection{Conclusions}

We associated an infinite-dimensional B\"acklund loop Lie algebra structure of 
the Ka\v{c}--Moody type with the Ishimori $(2+1)$-dimensional 
spin model in both the compact \cite{Ish84} and the noncompact versions \cite{LMS92}.
As a matter of principle, it would be interesting to work out homomorphisms of $\bar{\mathfrak{k}}$ with
quotient Lie algebras also in the case when $Y\neq 0, Z\neq 0, K\neq 0$. This subject is currently under
investigation. 

Notice that, generalizing to the 
$(2+1)$--dimensional case the ``inverse'' 
procedure (see \eg \cite{ALLPS94,ALLPS95,PaWi02}), by suitable {\em nonlinear} representations of
the quotient Ka\v{c}--Moody algebra
(\ref{KM}) in terms of pseudopotentials, one could
obtain a whole family of (2+1)-dimensional nonlinear spin field equations
containing the original Ishimori spin models (see \eg \cite{ALLPS94} for similar results concerning the
$(1+1)$--dimensional Heisenberg spin model). This turns out to be interesting within the study of `hidden'
symmetries of nonlinear
$\sigma$-models. Details will be the subject of a separate paper. Here we just stress that such a family should
not be confused {\em a priori} with the well known hierarchy associated {\em via} the Lax pair; in
fact, as it is well known, the latter comes from a {\em linear} representation of quotient Lie algebras.

\subsection*{Acknowledgments}
The author is indebted with Professors R.A. Leo and G. Soliani for having introduced her - as
the supervisors of her doctoral thesis at the University of Lecce - to this intriguing topic and for always
enlightening discussions. 


\end{document}